\documentclass[preprint,12pt]{elsarticle}
\usepackage{amsthm,amsmath,amssymb}
\usepackage{graphicx}
\usepackage[colorlinks=true,citecolor=black,linkcolor=black,urlcolor=blue]{hyperref}
\usepackage{mathrsfs}

\def\cvd{~\vbox{\hrule\hbox{%
     \vrule height1.3ex\hskip0.8ex\vrule}\hrule } }
\theoremstyle{plain}
\newtheorem{theorem}{Theorem}

\newtheorem{corollary}[theorem]{Corollary}

\theoremstyle{definition}

\theoremstyle{remark}
\newtheorem{remark}[theorem]{Remark}

\begin{document}

\title{Several extreme coefficients of the Tutte polynomial of graphs}
\author{Helin Gong\\
\small School of Mathematical Sciences\\[-0.8ex]
\small Xiamen University\\[-0.8ex]
\small P. R. China\\
\small Department of Fundamental Courses\\[-0.8ex]
\small Zhejiang Industry Polytechnic College\\[-0.8ex]
\small P. R. China\\
Mengchen Li\\
\small School of Mathematical Sciences\\[-0.8ex]
\small Xiamen University\\[-0.8ex]
\small P. R. China\\
Xian'an Jin\\
\small School of Mathematical Sciences\\[-0.8ex]
\small Xiamen University\\[-0.8ex]
\small P. R. China\\
\small\tt Email:xajin@xmu.edu.cn
}
\begin{abstract}
Let $t_{i,j}$ be the coefficient of $x^iy^j$ in the Tutte polynomial $T(G;x,y)$ of a connected bridgeless and loopless graph $G$ with order $n$ and size $m$. It is trivial that $t_{0,m-n+1}=1$ and $t_{n-1,0}=1$. In this paper, we obtain expressions of another eight extreme coefficients $t_{i,j}$'s with   $(i,j)=(0,m-n)$,$(0,m-n-1)$,$(n-2,0)$,$(n-3,0)$,$(1,m-n)$,$(1,m-n-1)$,$(n-2,1)$ and $(n-3,1)$ in terms of small substructures of $G$. Among them, the former four can be obtained by using coefficients of the highest, second highest and third highest terms of chromatic or flow polynomials, and vice versa. We also discuss their duality property and their specializations to extreme coefficients of the Jones polynomial.
\end{abstract}

\begin{keyword}
Extreme coefficients\sep Tutte polynomial\sep convolution formula\sep chromatic polynomial\sep flow polynomial \sep Jones polynomial

\MSC 05C31\sep 57M27
\end{keyword}
\maketitle

\section{Introduction}
\noindent

Throughout this paper we consider finite (undirected) graphs that allow parallel edges and may have loops.  Let $G = (V, E)$ be a graph with vertex set $V$ and edge set $E$. The order, size of $G$ and the number of connected components of $G$ are denoted by $n=n(G)$, $m=m(G)$ and $c=c(G)$, respectively. The complete graph, the empty graph, the path and the cycle of order $n$ is denoted by $K_n, E_n, P_n$ and $C_n$, respectively. For $A \subseteq E$, we denote by $G[A]$ the subgraph induced by $A$, $G/A$ the graph obtained from $G$ by contracting all edges in $A$,   $G-A$ the graph obtained from $G$ by deleting edges in $A$ and $G|_A $ the restriction of $G$ to $A$, namely $G|_A = G-(E\backslash A)$.

The Tutte polynomial $T(G; x, y)$ of a graph $G = (V,E)$, introduced in \cite{TUTTE}, is a two-variable polynomial which can be recursively defined as:
\begin{eqnarray*}\nonumber
T(G;x,y)=
\begin{cases}
1& \text{if $E = \emptyset$}\\
xT(G/e; x,y) & \text{if $e$ is a bridge}\\
yT(G-e; x,y) & \text{if $e$ is a loop}\\
T(G/e;x,y) + T(G-e;x,y)  & \text{if $e$ is neither a loop nor a bridge}.
\end{cases}
\end{eqnarray*}
It is independent of the order of edges selected for deletion and contraction in the reduction process to the empty graph. One way of seeing this is through the rank-nullity expansion of the Tutte polynomial. Let $A\subseteq E$. We identify $A$ with the spanning subgraph $(V, A)$ of $G$, i.e. $G|_A$,  temporarily for the sake of simplicity. Let $\rho(A)$ denote the rank $n-c(A)$, $\gamma(A)$ denote the nullity $|A|-n+c(A)$. Then
\begin{eqnarray}\nonumber
T (G; x, y) = \sum_{A \subseteq E}(x- 1)^{\rho(E)-\rho(A)}(y-1)^{\gamma(A)}.
\end{eqnarray}
Moreover, the Tutte polynomial has a spanning forest expansion \cite{TUTTE}, i.e.
\begin{eqnarray*}
T (G; x, y) = \sum_{i,j}t_{ij}x^iy^j,
\end{eqnarray*}
where $t_{ij}$ is the number of spanning forests of $G$  with internal activity $i$ and external activity $j$.

The Tutte polynomial contains as special cases the chromatic polynomial $P(G; \lambda)$ which counts proper $\lambda$-colorings of $G$ and the flow polynomial $F(G; \lambda)$ which counts nowhere-zero $AG$-flows of $G$, where $AG$ is a finite Abelian group and $|AG|=\lambda$. Namely,
\begin{eqnarray}
P(G; \lambda)&=&(-1)^{\rho(E)}\lambda^{c}T(G; 1 - \lambda, 0),\label{chro}\\
F(G; \lambda)&=&(-1)^{\gamma(E)}T(G; 0, 1-\lambda).\label{flow}
\end{eqnarray}
In \cite{KOOK}, Kook obtained the following convolution formula for the Tutte polynomial, which will be used in Section 2.
\begin{eqnarray}
T(G;x,y) = \sum_{A\subseteq E}T(G/A; x,0)T(G|_A;0,y).\label{con}
\end{eqnarray}

It is obvious that $t_{00} = 0$ if $|E|>0$. It is proved that $t_{01} = t_{10}$ if $|E|>1$ \cite{BOLLOBAS}, and it is called $\beta$ invariant. $\beta \neq 0$ implies that the considered graph is loopless and 2-connected \cite{BRYLAWSKI} and the $\beta$ invariant enumerates the bipolar orientations of a rooted 2-connected plane graph \cite{FOR}. For surveys of results and applications of the Tutte polynomial, we refer the readers to \cite{BRYLAWSKI,WELSH,ELLIS}. It is basic in graph theory to establish relations between the coefficients of graph polynomials and subgraph structures in the graph. See, for example, \cite{Biggs} for results on characteristic polynomial and chromatic polynomial. The purpose of the paper is to establish a relation between several extreme coefficients of the Tutte polynomial and subgraph structures of the graph. Let $G=(V,E)$ be a connected bridgeless and loopless graph. To state our results, we need some additional definitions and notations.

A \emph{parallel class} of $G$ is a
maximal subset of $E$ sharing the same endvertices. A parallel class is called \emph{trivial} if it contains only one edge. It is obvious that parallel classes partition the edge set $E$. A \emph{series class} of $G$ is a maximal subset $C$ of $E$ such
that the removal of any two edges from $C$ will increase the number of
connected components of the graph. Let $C\subseteq E$. Then $C$ is a series class if and only if $c(G-C)=|C|$ and $G-C$ is bridgeless. A series class is called \emph{trivial} if it contains only one edge. If a bridgeless and loopless graph $G$ is disconnected, then its series classes are defined to be the union of series classes of connected components of $G$. Series classes also partition the edge set $E$. See \cite{DJIN} for details.

Let $k_1,k_2,k_3\geq 1$ be integers. We denote by $\Theta_{k_1,k_2,k_3}$, the graph with two vertices $u$, $v$ connected by three internally disjoint paths of lengths $k_1,k_2$ and $k_3$. If $C\subseteq E$ satisfies: (1) $C=C_1\cup C_2\cup C_3$, $C_{i}\subseteq E$ and $|C_i|=k_i$ $(i=1,2,3)$, (2) $G$ has the structure as shown in Figure 1, and (3) $G-C= G_1\cup G_2\cup\cdots\cup G_{k_1+k_2+k_3-1}$ and each $G_i$ is connected and bridgeless for $i=1,2,\cdots,k_1+k_2+k_3-1$, then we say $C$ is a $\Theta$ class of $G$.  The total number of $\Theta$ classes of $G$ is denoted by $\Theta(G)$.

Let $k_1,k_2\geq 1$ be integers. We denote by $\infty_{k_1,k_2}$ the graph formed from two cycles $C_{k_1}$ and $C_{k_2}$ by identifying one vertex of $C_{k_1}$ with one vertex of $C_{k_2}$. If $C\subseteq E$ satisfies: (1) $C=C_1\cup C_2$, $C_{i}\subseteq E$ and $|C_i|=k_i$ $(i=1,2)$, (2) $G$ has the structure as shown in Figure 2, and (3) $G-C= G_1\cup G_2\cup\cdots\cup G_{k_1+k_2-1}$ and each $G_i$ is connected and bridgeless for $i=1,2,\cdots,k_1+k_2-1$, then we say $C$ is an $\infty$ class of $G$. The total number of $\infty$ classes of $G$ is denoted by $\infty(G)$.

\begin{figure}[htbp]
\label{fig1}
\centering
\includegraphics[width=.5\textwidth]{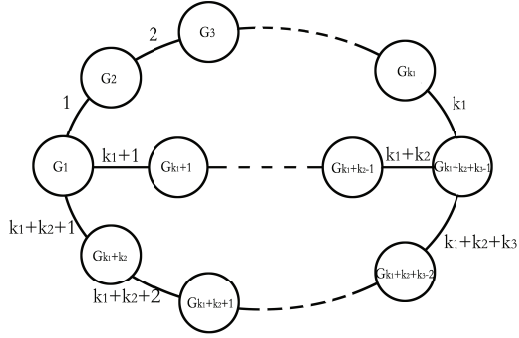}
\caption{A $\Theta$ class.}
\end{figure}

\begin{figure}[htbp]
\label{fig2}
\centering
\includegraphics[width=.5\textwidth]{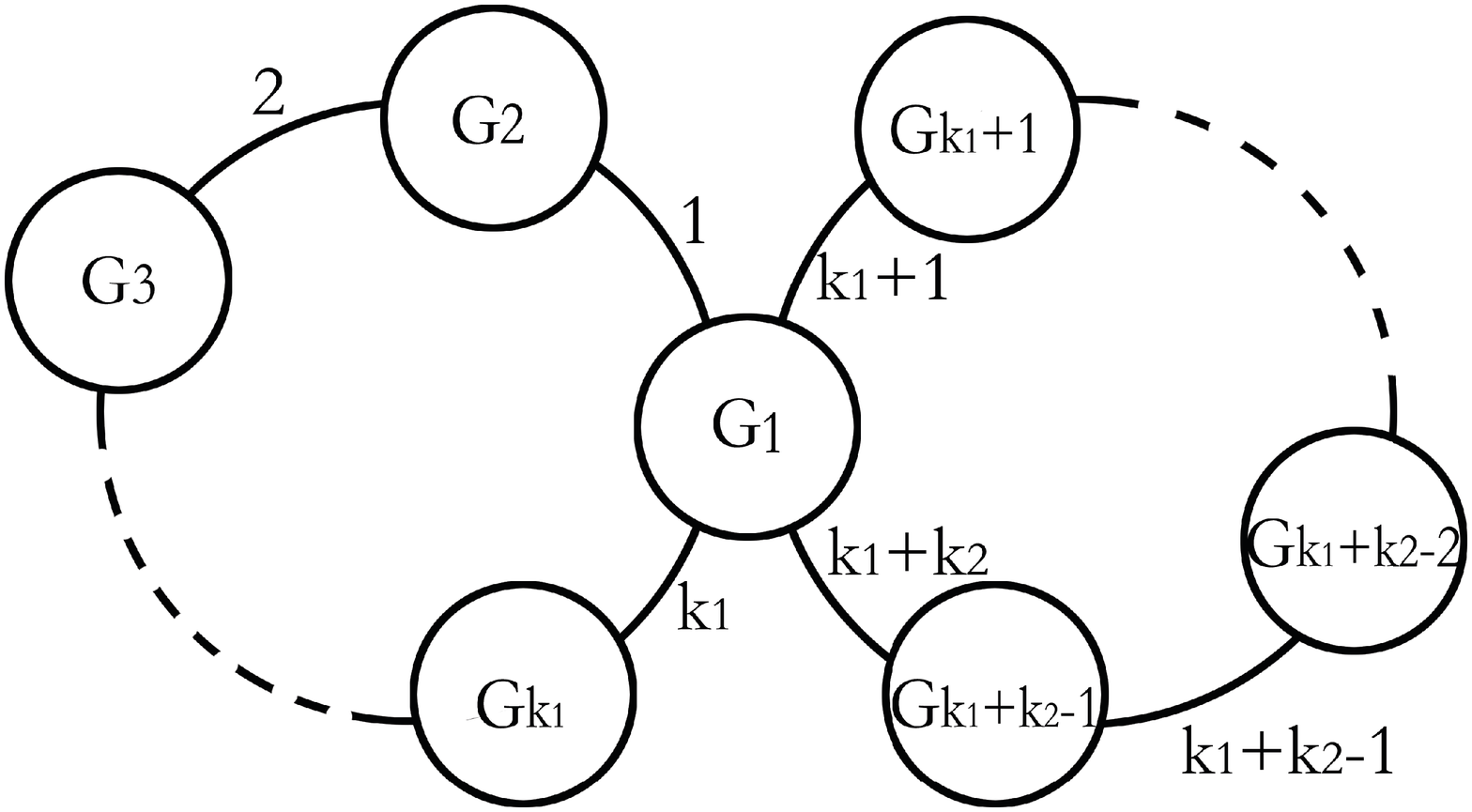}\label{eg}
\caption{An $\infty$ class.}
\end{figure}

Now we are in a position to state our main results.

\begin{theorem}\label{mainI}
Let $G = (V, E)$ be a loopless and bridgeless connected graph. Let $s(G)$ and  $s^*(G)$ be the number of series classes and non-trivial series classes of $G$, respectively.  Then
\begin{itemize}
\item[(1)] $t_{0,m-n+1}= 1$,
\item[(2)] $t_{0,m-n} = n + s(G) - m -1$,
\item[(3)] $t_{0,m-n-1} ={m-n+1 \choose 2} - (m-n)s(G) + {s(G) \choose 2} - \Theta(G)$,
\item[(4)] $t_{1,m-n} =s^*(G)$, \textrm{and}
\item[(5)]
$
t_{1,m-n-1} =-s^*(G)(m-n)+\sum_{\stackrel{A\subseteq E}{A\ \text{is a nontrivial series class} }  }s(G-A) + \Theta(G).
$
\end{itemize}
\end{theorem}

\begin{theorem}\label{mainII}
Let $G = (V, E)$ be a loopless and bridgeless connected graph. Let $p(G)$ and $p^*(G)$ be the number of parallel classes and non-trivial parallel classes of $G$, respectively. Let $\Delta(\tilde{G})$ be the number of triangles of $\tilde{G}$, the graph obtained from $G$ by replacing each parallel class by a single edge. Then
\begin{itemize}
\item[(1)] $t_{n-1,0} = 1$,
\item[(2)] $t_{n-2,0} =  p(G) - n +1$,
\item[(3)] $t_{n-3,0} ={n-1 \choose 2} - (n-2)p(G) + {p(G) \choose 2}-\Delta(\tilde{G})$,
\item[(4)] $t_{n-2,1} =p^*(G)$, \textrm{and}
\item[(5)] $t_{n-3,1} = -p^*(G)(n-2)+\sum_{\stackrel{A\subseteq E}{ A\ \text{is a nontrivial  parallel class} }  }p(G/A) + \Delta(\tilde{G}).$
\end{itemize}
\end{theorem}

\begin{remark}
$t_{0,m-n+1}, t_{0,m-n}$ and $t_{0,m-n-1}$ can be obtained from coefficients of the flow polynomial $F(G; 1-y)$ and $t_{n-1,0}, t_{n-2,0}$ and $t_{n-3,0}$ can be obtained from coefficients of the chromatic polynomial $P(G;1-x)$. These will be seen from proofs of Theorems 1 and 2 and further clarified in Section 3. Another four coefficients are completely new as far as we know.
\end{remark}

\section{Proofs}
\noindent

To prove Theorems \ref{mainI} and \ref{mainII}, we need two results that express the chromatic and flow polynomials of graphs as characteristic polynomials of some lattice related to graphs, respectively, which was proven by G.~C. Rota in the 1960s \cite{ROTA}.

\subsection{M\"{o}bius function, the chromatic and flow polynomials}
\noindent

Let $P$ be a poset (a finite set $S$ partially ordered by the relation $\leq$). The unique \emph{minimum element} and unique \emph{maximum element} in $P$, if they exist, are denoted by $\widehat{0}=\widehat{0}_P, \widehat{1}=\widehat{1}_P$, respectively. A segment $[x,y]$, for $x,y\in P$, is the set of all elements $z$ between $x$ and $y$, i.e. $\{z | x\leq z\leq y\}$. Note that the segment $[x,y]$ endowed with the induced order structure is a poset in its own right and $\widehat{0}_{[x,y]} = x, \widehat{1}_{[x,y]} = y$. An element $y$ \emph{covers} an element $x$ when the segment $[x,y]$ contains two elements. A poset is \emph{locally finite} if every segment is finite. Let $P$ be a locally finite poset. Then
the M\"{o}bius function of $P$ is an integer-valued function defined on the Cartesian set $P\times P$ such that
\begin{align*}
\mu(x,y) =1\ & \mbox{if}\ x =  y\\
\mu(x,y) =0\ & \mbox{if}\ x \nleq  y\\
\sum_{x\leq z\leq y} \mu(x,z)=0\ & \mbox{if}\  x <y.
\end{align*}
If $P$ has a $\widehat{0}$, then
\begin{eqnarray}
\mu(\widehat{0},x) = -\sum_{y<x} \mu(\widehat{0},y)\ \mbox{if}\  x >\widehat{0}.
\end{eqnarray}

A finite poset $P$ is \emph{ranked (or graded)} if for every $x\in P$ every maximal chain
with $x$ as top element has the same length, denoted $rk(x)$. Here the length of
a chain with $k$ elements is $k - 1$. If $P$ is ranked, the function $rk$ called the
\emph{rank function}, is zero for minimal elements of $P$ and $rk(y) = rk(x) + 1$ if
$x, y\in P$ and $y$ covers $x$.

Let $P$ be a ranked poset and $t$ be a variable. Then the characteristic polynomial of $P$ is defined by
\begin{eqnarray}
q(P;t) = \sum_{x\in P}\mu(\widehat{0},x)t^{rk(\widehat{1})-rk(x)}.
\end{eqnarray}

Let $G = (V, E)$ be a loopless graph.  A \emph{bond} of $G$ is a spanning subgraph $H \subseteq G$ such that each
connected component of $H$ is a vertex-induced subgraph of $G$. Then the set $L(G)$ consisting of all bonds of $G$ forms a graded lattice ordered by the refinement relation on the set of partitions of $V$, that is,  $K \in L(G) \leq H\in L(G)$  means that $\{V(K_1), \cdots, V(K_s)\}$ is finer than $\{V(H_1), \cdots, V(H_t)\} $, where $K_1, \cdots, K_s$ are connected components of $K$ and $H_1, \cdots, H_t$ are connected components of $H$. Moreover, $rk(H) = \rho(H)$ for $H\in L(G)$ and $P(G; \lambda ) = \lambda^{c(G)}q(L(G); \lambda ) $.

\begin{theorem}[\cite{ROTA}] Let $G$ be a loopless graph. Then
\label{thA}
\begin{align}
P(G; \lambda )=  \sum_{H\in L(G)}\mu(E_n,H)\lambda^{c(H)}.\label{ty}
\end{align}
\end{theorem}

\begin{remark}
Note that when $G$ contains parallel edges, Theorem \ref{thA} is still valid if we take $\tilde{G}$ in place of $G$ in the right side of (\ref{ty}).
\end{remark}

Let $G = (V,E)$ be a bridgeless graph. The set
$L'(G)$ consisting of all spanning subgraphs of $G$ without bridges also forms a graded lattice with the partial order defined by $H_1 \leq  H_2$ if $E(H_1) \supseteq E(H_2)$. Moreover, $rk(H) = \gamma(G) - \gamma(H)$ for $H\in L'(G)$.

\begin{theorem} Let $G$ be a bridgeless graph. Then
\label{thB}
\begin{align}
F(G; \lambda) = \sum_{H\in L'(G)}\mu(G,H)\lambda^{\gamma(H)}.
\end{align}
\end{theorem}
\noindent\textbf{Proof.} Let $H \in L'(G)$.  Suppose $N_{=}(H)$ is the function counting $AG$-flows of $\overrightarrow{G}$ such that $\mathbf{0}$ is assigned exactly on edges of  $E\backslash E(H)$, and $N_{\geq}(H) =\sum_{H'\geq H}N_=(H')$ is the function counting $AG$-flows of $\overrightarrow{G}$ such that  $\mathbf{0}$ is assigned  at least on edges of $E\backslash E(H)$. Note that $N_=(G) = F(G; \lambda)$.  By the M\"{o}bius Inversion Theorem \cite{BENDER},

\[
N_{=}(G) = \sum_{H\in L'(G)}\mu(G,H)N_{\geq}(H).
\]
It is not difficult to see that $N_{\geq}(H) = \lambda^{\gamma(H)}$, which completes the proof. \hfill\cvd

\vskip0.2cm

We write $\omega_1 = 1-x$ and $\omega_2 = 1-y$. Keep in mind that $\widehat{0}$ will be the graph $G$ itself when $L'(G)$ is concerned while $\widehat{0}$ will be the empty graph of order $n(G)$ when $L(G)$ is concerned. By inserting Eqs. (\ref{chro}) and (\ref{flow}) into Eq. (\ref{con}) , we obtain
\begin{eqnarray}
&&T(G;x,y)=\sum_{A\subseteq E}T(G/A; x,0)T(G|_A;0,y)\nonumber\\
&=& \sum_{A\subseteq E}(-1)^{n(G/A)-c(G/A)+|A|-n+c(G|_A)}\omega_1^{-c(G/A)}P(G/A;\omega_1)F(G|_A;\omega_2).
\label{eq-sum}
\end{eqnarray}
Thus, we only consider $A$'s in the RHS of (\ref{eq-sum}) such that $G/A$ is loopless and $G|_A$ is bridgeless in the proofs of Theorems \ref{mainI} and \ref{mainII}. Otherwise, $P(G/A;\omega_1)=0$ or $F(G|_A;\omega_2)=0$.

The following notations will be used in the next two subsections. Let $D_k$ be the dipole graph of size $k$, i.e.  two distinct vertices connected by $k$ parallel edges. Let $P_{k_1, k_2}$~($k_i\geq 1$ for each $i=1,2$) be the \emph{multi-path} with $3$ vertices $v_1,v_2,v_{3}$ connected by $k_i$ parallel edges between $v_i$ and $v_{i+1}$ for $i=1,2$.
Let $C_{k_1, k_2, k_3}$ ($k_i \geq 1$ for each $i=1,2,3$) be the \emph{multi-cycle} with $3$ vertices $v_1, v_2, v_{3}$ connected by $k_i$ parallel edges between $v_i$ and $v_{i+1}$ for $i=1,2$ and $k_3$ parallel edges between $v_3$ and $v_1$.

%%%%%%%%%%%%%%%%%%%%%%%%%%%%%%%%%%%%%%%%%%%
% Theorem 1.1
%%%%%%%%%%%%%%%%%%%%%%%%%%%%%%%%%%%%%%%%%%%
\subsection{Proof for Theorem \ref{mainI}}
\noindent

Note that the degree in $y$ in Theorem \ref{mainI} is $m-n+1, m-n$ and $m-n-1$. We shall consider $A$'s with $G|_A$ close to $G$.

\noindent\textbf{Case 1} $ |A| = m$.  In this case, $G/A=K_1$ and $ G|_A=G$. The corresponding contribution of $A$ to the summation of (\ref{eq-sum}) is
\begin{eqnarray}
(-1)^{m-n+1}[\omega_2^{m-n+1} + \sum_{\stackrel {H\in L'(G)}{rk(H)=1}}\mu(\widehat{0},H)\omega_2^{m-n}+\sum_{\stackrel {H\in L'(G)}{ rk(H)=2}}\mu(\widehat{0},H)\omega_2^{m-n-1}+\cdots].
\label{eq-A1}
\end{eqnarray}

\noindent\textbf{Case 2} $|A| = m-1$. In this case, $ G/A $ will contain a loop.

\noindent\textbf{Case 3} $|A| = m-2$. In this case, $E\backslash A$ is exactly a series class with cardinality $2$ of $G$. Moreover, $ G/A=D_2$ and hence $P(G/A; \omega_1)= -x(1-x)$. It is clear that $c(G/A) = 1$ and $c(G|_A) =  c(G-(E\backslash A)) = 2$. Hence,
\[(-1)^{n(G/A)-c(G/A)+|A|-n+c(G|_A)}\omega_1^{-c(G/A)} = (-1)^{m-n+1}(1-x)^{-1}. \]
Thus, the contribution of $A$ can be written as
\begin{align}
(-1)^{m-n}x[\omega_2^{m-n} + \sum_{\stackrel {H\in L'(G|_A)}{ rk(H)=1}}\mu(\widehat{0},H)\omega_2^{m-n-1}+\cdots ]
\label{eq-A2}
\end{align}
since $\deg(F(G|_A; \lambda))= \gamma(G|_A) = m - n$.

\noindent\textbf{Case 4} $|A| = m-3$. There are two subcases.

\noindent\textbf{(a)} $c(G|_A) = c(G - E\backslash A) = 3$. Then $E\backslash A$ is exactly a series class with cardinality 3. Moreover, $ G/A =  C_3$. Then $n(G/A) =3, c(G/A) = 1$, $c(G|_A) = 3$ and $P(G/A; \omega_1)= x(1-x)(1+x)$. So
\[(-1)^{n(G/A)-c(G/A)+|A|-n+c(G|_A)}\omega_1^{-c(G/A)} = (-1)^{m-n}(1-x)^{-1}. \]
Thus, the contribution of $A$ can be written as
\begin{align}
(-1)^{m-n}x(1+x)[\omega_2^{m-n} + \sum_{\stackrel {H\in L'(G|_A)}{ rk(H)=1}}\mu(\widehat{0},H)\omega_2^{m-n-1}+\cdots ].
\label{eq-A4}
\end{align}
since $\deg(F(G|_A; \lambda))= \gamma(G|_A) = m - n$.

\noindent\textbf{(b)} $c(G|_A) = c(G - E\backslash A) = 2$. Then $E\backslash A$ a $\Theta$ class. Moreover, $G/A = \Theta_{1,1,1}$. Note that $n(G/A) =2$, $c(G/A) = 1$, $c(G|_A) = 2$ and $P(G/A; \omega_1)= -x(1-x)$. So
\[(-1)^{n(G/A)-c(G/A)+|A|-n+c(G|_A)}\omega_1^{-c(G/A)} = (-1)^{m-n}(1-x)^{-1}. \]
Thus, the contribution of $A$ can be written as
\begin{align}
(-1)^{m-n-1}x[\omega_2^{m-n-1} +\cdots ]
\label{eq-A3}
\end{align}
since $\deg(F(G|_A; \lambda))= \gamma(G|_A) = m - n-1$.

\noindent\textbf{Case 5}  $|A| = m- k$ ($k\geq 4$). Since the degree in $y$ we considered in Theorem \ref{mainI} is at least $m-n-1$, only the following two subcases are needed to be considered.

\noindent\textbf{(a)} $c(G - E\backslash A) = k$. It means that $E\backslash A$ is exactly a series class with cardinality $k$. Note that $G/A=C_k$. Then the contribution of $A$ can be written as
\begin{align}
&(-1)^{k-1+m-k-n+ k}\omega_1^{-1}P(C_k;\omega_1)[\omega_2^{m - n}+ \sum_{\stackrel {H\in L'(G|_A)}{ rk(H)=1}}\mu(\widehat{0},H)\omega_2^{m-n-1}+\cdots]\nonumber\\
=&(-1)^{m-n+k-1}(1-x)^{-1}[(-1)^k(-x) + (-x)^k][\omega_2^{m-n}+\sum_{\stackrel {H\in L'(G|_A)}{ rk(H)=1}}\mu(\widehat{0},H)\omega_2^{m-n-1}+\cdots] \nonumber\\
=&(-1)^{m-n}[x+\cdots+x^{k-1}][\omega_2^{m-n}+\sum_{\stackrel {H\in L'(G|_A)}{ rk(H)=1}}\mu(\widehat{0},H)\omega_2^{m-n-1}+\cdots] \nonumber\\
=& xy^{m-n} - [(m-n)+ \sum_{\stackrel {H\in L'(G|_A)}{ rk(H)=1}}\mu(\widehat{0},H)]xy^{m-n-1}+\cdots.\label{eq-A5}
\end{align}

\noindent\textbf{(b)} $c(G - E\backslash A) = k-1$. Then $E\backslash A$ is either a $\Theta$ class or a $\infty$ class.

\noindent\textbf{(1)} For the first case, $G/A$ will be the theta graph $\Theta_{k_1,k_2,k_2}$ with $k_1, k_2, k_3\geq 1$ and $k_1 +k_2 +k_3 = k$. Then we know from \cite{BH} that
\begin{align*}
P(G/A;\omega_1)
=&\frac{\operatorname*{\prod}\limits_{i=1}^3[(\omega_1-1)^{k_i+1}+(-1)^{k_i+1}(\omega_1-1)]}{[\omega_1(\omega_1-1)]^2}\nonumber\\
&+\frac{\operatorname*{\prod}\limits_{i=1}^3[(\omega_1-1)^{k_i}+(-1)^{k_i}(\omega_1-1)]}{\omega_1^2}.
\end{align*}

If $k_1, k_2, k_3\geq 2$, the chromatic polynomial of $G/A$ can be written as
\begin{align*}
P(G/A;\omega_1)=&(-1)^kx(1-x)(1+\cdot\cdot\cdot+x^{k_1-1})(1+\cdot\cdot\cdot+x^{k_2-1})(1+\cdot\cdot\cdot+x^{k_3-1})\nonumber\\
&+(-1)^{k+1}x^3(1-x)(1+\cdot\cdot\cdot+x^{k_1-2})(1+\cdot\cdot\cdot+x^{k_2-2})(1+\cdot\cdot\cdot+x^{k_3-2}).
\end{align*}

Then the contribution of $A$ to the summation of (\ref{eq-sum}) is
\begin{align*}
&(-1)^{(k-1)-1+m-k-n+(k-1)}\omega_1^{-1}P(G/A;\omega_1)[\omega_2^{m-n-1}+\cdots]\nonumber\\
=&(-1)^{m-n-1}x(1+\cdots+x^{k_1-1})(1+\cdots+x^{k_2-1})(1+\cdots+x^{k_3-1})[\omega_2^{m-n-1}+\cdots]\nonumber\\
&+(-1)^{m-n}x^3(1+\cdots+x^{k_1-2})(1+\cdots+x^{k_2-2})(1+\cdots+x^{k_3-2})[\omega_2^{m-n-1}+\cdots].
\end{align*}

If there is one $k_i$ such that $k_i=1$, then
\begin{align*}
P(G/A;\omega_1)=(-1)^kx(1-x)(1+\cdot\cdot\cdot+x^{k_1-1})(1+\cdot\cdot\cdot+x^{k_2-1})(1+\cdot\cdot\cdot+x^{k_3-1}).
\end{align*}

The contribution of $A$ is
\begin{align*}
&(-1)^{(k-1)-1+m-k-n+(k-1)}\omega_1^{-1}P(G/A;\omega_1)[\omega_2^{m-n-1}+\cdots]\nonumber\\
=&(-1)^{m-n-1}x[1+\cdots+x^{k_1-1}][1+\cdots+x^{k_2-1}][1+\cdots+x^{k_3-1}][\omega_2^{m - n-1}+\cdots].
\end{align*}

It can be seen that no mater what values of $k_1$, $k_2$, $k_3$ are, the contributions of $A$ are the same, i.e.
\begin{align}
(-1)^{m-n-1}x[1+\cdots+x^{k_1-1}][1+\cdots+x^{k_2-1}][1+\cdots+x^{k_3-1}][\omega_2^{m - n-1}+\cdots]\label{eq-A7}.
\end{align}

Note that Eq. (\ref{eq-A7}) coincides with Eq. (\ref{eq-A3}).

\noindent\textbf{(2)} For the second case, $G/A$ will be $\infty_{k_1,k_2}$ with $k_1,k_2\geq 1$, $k_1 +k_2 = k$. If $k_1$ or $k_2$ is 1, then $G/A$ contains a loop. So we only consider $k_1,k_2\geq 2$. Then the contribution of $A$ can be written as
\begin{align*}
&(-1)^{(k-1)-1+m-k-n+(k-1)}\omega_1^{-1}\frac{C(C_{k_1};\omega_1)C(C_{k_2};\omega_1)}{\omega_1}[\omega_2^{m-n-1}+\cdot\cdot\cdot]\nonumber \\
=&(-1)^{m-n-1}x^2[1+\cdot\cdot\cdot+x^{k_1-2}][1+\cdot\cdot\cdot+x^{k_2-2}][\omega_2^{m - n-1}+\cdots].
\end{align*}

Now we combine all above cases to obtain Theorem \ref{mainI}.
Note that we only need to consider Case 1 to determine $t_{0,m-n+1},t_{0,m-n}$ and $t_{0,m-n-1}$ since Case 2 has no contribution and contributions of Cases 3 to 5 and their subcases and subsubcases all include the variable $x$. To determine $t_{1,m-n}$, we only need consider Case 3, Case 4(a) and Case 5(a). To determine $t_{1,m-n-1}$, we need consider the first type: Case 3, Case 4(a) and Case 5(a) and the second type: Case 4(b) and Case 5(b)(1). Note that Case 5(b)(2) includes $x^2$ and has no contribution. Now let's insert (\ref{eq-A1})-(\ref{eq-A7}) into (\ref{eq-sum}). We obtain:
\begin{align*}
&T(G; x, y)\nonumber\\
=&y^{m-n+1} + (-1)[(m-n+1) + \sum_{\stackrel {H\in L'(G)}{rk(H)=1}}\mu(\widehat{0},H)]y^{m-n}  \\
& + [{m-n+1 \choose 2} + (m-n)\sum_{\stackrel {H\in L'(G)}{rk(H)=1}}\mu(\widehat{0},H) + \sum_{\stackrel {H\in L'(G)}{rk(H)=2}}\mu(\widehat{0},H)]y^{m-n-1}\\
& + \cdots \\
&+ \sum_{\stackrel {A \subseteq E }{E\backslash A\ \text{is a nontrivial serial class}}} xy^{m-n}\\
&+ \Big\{\sum_{\stackrel {A \subseteq E }{E\backslash A\  \text{is a nontrivial serial class}}}[-(m-n)-\sum_{\stackrel {H\in L'(G|_A)}{rk(H)=1}}\mu(\widehat{0},H)]+ \Theta(G)\Big\}xy^{m-n-1}\\
& + \cdots
\end{align*}

Let $H=G-A \in L'(G)$. Then $rk(H) =1$ $\Longleftrightarrow$ $\gamma(H)=m-n$ $\Longleftrightarrow$ $|A| = c(G -A)$. Thus
\[\sum_{\stackrel {H\in L'(G)}{rk(H)=1}} \mu(\widehat{0}, H)=\sum_{\stackrel {A\subseteq E, |A| = c(G-A)} {G - A\  \text{is bridgeless}}} (-1)= \sum_{A\  \text{is a series class}} (-1)=-s(G). \]

\begin{figure}[htbp]
\label{fig2}
\centering
\includegraphics[width=.5\textwidth]{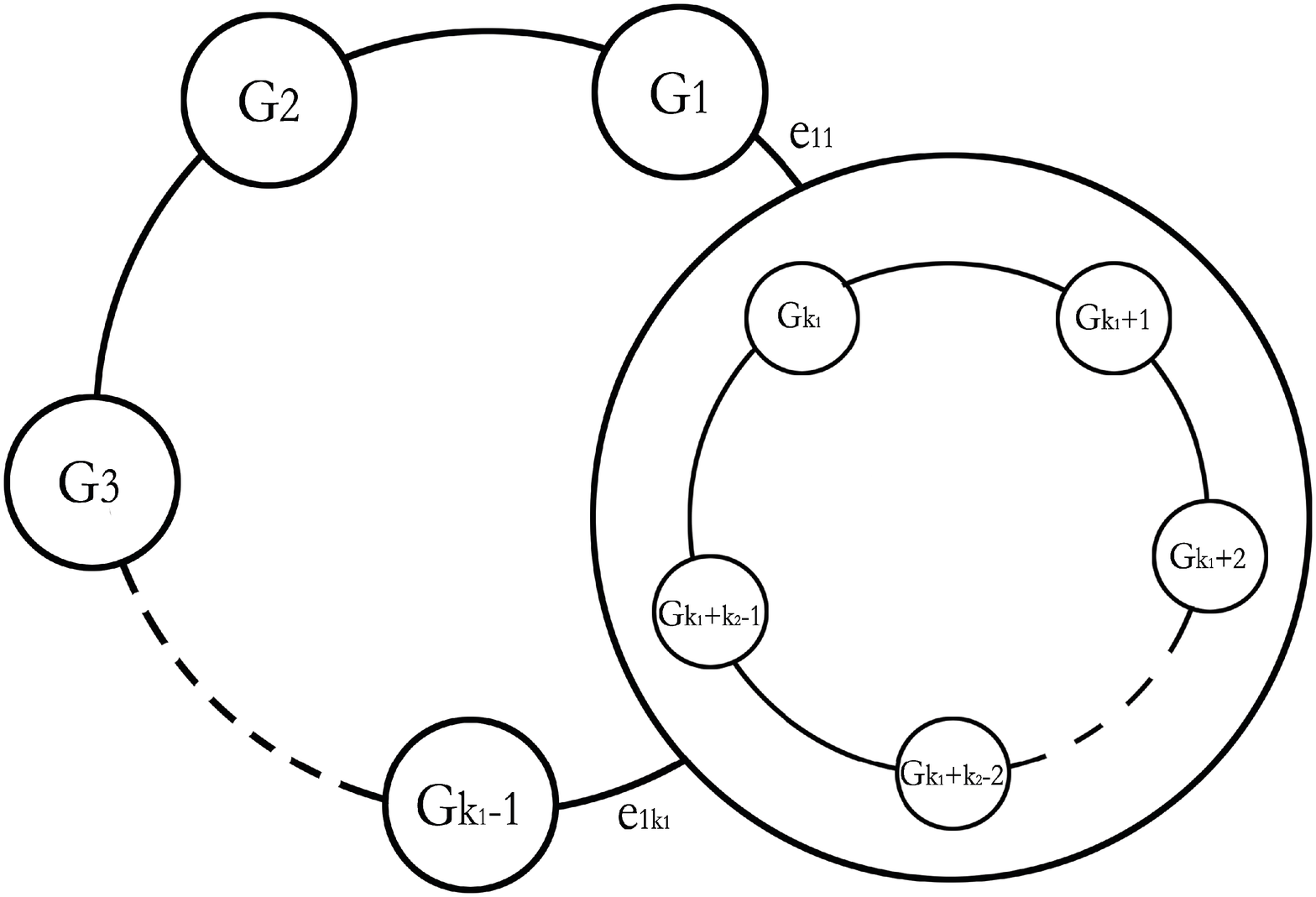}
\caption{$H'=G_1\cup\cdots\cup G_{k_1-1}\cup G_{k_1}\cup\cdots\cup G_{k_1+k_2-1}$.}
\end{figure}

Similarly, let $H'=G-A' \in L'(G)$. Then $rk(H') =2$ $\Longleftrightarrow$ $|A'| = c(G -A') +1 $ $\Longleftrightarrow$ $A'-A$ is a series class of $G-A$ $\Longleftrightarrow$ $G$ has the structure like Figure 3 $\Longleftrightarrow$ $A'$ is either a $\Theta$ class or a $\infty$ class. If $A'$ is a $\Theta$ class, then $\mu(\hat{0},H')=2$ and if $A'$ is an $\infty$ class, then $\mu(\hat{0},H')=1$. Note that $\infty(G)=\binom{s(G)}{2}-3\Theta(G)$. Thus
\begin{align*}
\sum\limits_{\scriptstyle H\in L'(G) \atop \scriptstyle rk(H)=2}\mu(\hat{0},H)=1\times\left[\binom{s(G)}{2}-3\Theta(G)\right]+2\times\Theta(G)=\binom{s(G)}{2}-\Theta(G).
\end{align*}
This completes the proof of Theorem \ref{mainI}.
\hfill\cvd

\vskip0.2cm
%%%%%%%%%%%%%%%%%%%%%%%%%%%%%%%%%%%%%%%%%%%
% Theorem 1.2
%%%%%%%%%%%%%%%%%%%%%%%%%%%%%%%%%%%%%%%%%%%

\subsection{Proof for Theorem \ref{mainII}}
\noindent

We only give a sketch of the proof of Theorem \ref{mainII}.

\noindent\textbf{Case 1}  $ A = \emptyset $. In this case, $ G/A =G$ and $ G|_A =E_n$. The contribution of $A$ to the summation (\ref{eq-sum}) is
\begin{eqnarray}
(-1)^{n-1}[\omega_1^{n-1} + \sum_{\stackrel {H\in L(G)}{rk(H)=1}}\mu(\widehat{0},H)\omega_1^{n-2}+\sum_{\stackrel {H\in L(G)}{ rk(H)=2}}\mu(\widehat{0},H)\omega_1^{n-3}+\cdots].\label{ccc}
\label{eq-B1}
\end{eqnarray}

\noindent\textbf{Case 2} $|A| = 1$. $G|_A$ will have a bridge.

\noindent\textbf{Case 3} $|A| = 2$. $A$ is exactly a parallel class of cardinality 2, $G|_A =D_2\cup E_{n-2}$ and
\[(-1)^{n(G/A)-c(G/A)+|A|-n+c(G|_A)}\omega_1^{-c(G/A)} = (-1)^{n-1}\omega_1^{-1}. \]
It follows that the contribution of $A$ can be written as
\begin{align}
(-1)^{n-1}(-y)[\omega_1^{n-2} + \sum_{\stackrel {H\in L(G/A)}{ rk(H)=1}}\mu(\widehat{0},H)\omega_1^{n-3}+\cdots ].
\label{eq-B2}
\end{align}

\noindent\textbf{Case 4} $|A| = 3$. There are only two subcases to be considered.

\noindent\textbf{(a)} $A$ is a parallel class. In this case $n(G/A) =n-1$, $c(G/A) = 1$, $c(G|_A) = n-1$. Hence
\[(-1)^{n(G/A)-c(G/A)+|A|-n+c(G|_A)}\omega_1^{-c(G/A)} = (-1)^{n}\omega_1^{-1}. \]
Thus, the contribution of $A$ can be written as
\begin{align}
(-1)^{n}(y^2+y)[\omega_1^{n-2} + \sum_{\stackrel {H\in L(G/A)}{ rk(H)=1}}\mu(\widehat{0},H)\omega_1^{n-3}+\cdots ].
\label{eq-B3}
\end{align}
since $\deg(P(G/A;\lambda))= n-1$ and $F(G|_A; \omega_2) = (1-\omega_2) + (1-\omega_2)^2 = y^2 +y $.

\noindent\textbf{(b)} Each edge of $A$ is a trivial parallel class and $G[A]=C_3$. In this case,  $n(G/A) =n-2$, $c(G/A) = 1$, $c(G|_A) = n-2$. Thus,
\[(-1)^{n(G/A)-c(G/A)+|A|-n+c(G|_A)}\omega_1^{-c(G/A)} = (-1)^{n}\omega_1^{-1}. \]
Note that $\deg(P(G/A; \lambda))= n-2$ and $F(G|_A; \omega_2) = -y$, the contribution of $A$ can be written as
\begin{align}
(-1)^{n}(-y)[\omega_1^{n-3} + \cdots ].
\label{eq-B4}
\end{align}

\noindent\textbf{Case 5} $|A| = k$ ($k\geq 4$). Since the degree of $x$ considered in Theorem \ref{mainII} is at least $n-3$,  we only need to discuss the following subcases.

\noindent\textbf{(a)} $A$ is a  parallel class. Then the contribution of $A$ is
\begin{align}
&(-1)^{n(G/A)-c(G/A)+|A|-n+c(G|_A)}\omega_1^{-c(G/A)}F(G|_A;\omega_2)P(G/A; \omega_1)\nonumber\\
=&(-1)^{(n-1)-1+k-n+(n-1)}F(D_k;\omega_2)[\omega_1^{n-2}+\sum_{\stackrel {H\in L(G/A)}{ rk(H)=1}}\mu(\widehat{0},H)\omega_1^{n-3}+\cdots]\nonumber\\
=&(-1)^{n+k-1}[(-1)^{k-1}\sum_{i=1}^{k-1}y^{i}][\omega_1^{n-2}+\sum_{\stackrel {H\in L(G/A)}{ rk(H)=1}}\mu(\widehat{0},H)\omega_1^{n-3}+\cdots]\nonumber\\
=& yx^{n-2} - [(n-2)+ \sum_{\stackrel {H\in L(G/A)}{ rk(H)=1}}\mu(\widehat{0},H)] yx^{n-3} +\cdots.\label{eq-B5}
\end{align}

\noindent\textbf{(b)} $G[A]=P_{k_1, k_2}$ or $D_{k_1}\cup D_{k_2}$ ($k_i\geq 2, k_1+k_2 = k$) and $G/A$ has no loops. Since the lowest degree in $y$ of the polynomial $F(G|_A;\omega_2)$ is 2. So we need not consider such $A$'s.

\noindent\textbf{(c)} $G[A]=C_{k_1,k_2,k_3}$ ($k_i\geq 1, k_1+k_2+k_3=k$) and $G/A$ has no loops. Then the contribution of $A$ is
\begin{align}
&(-1)^{n(G/A)-c(G/A)+|A|-n+c(G|_A)}F(C_{k_1, k_2, k_3}; \omega_2)\omega_1^{-1}[\omega_1^{n-2}+\cdots]\nonumber\\
=&(-1)^{n-k+1}[(-1)^{k}(y + 3y^2 + \cdots)][\omega_1^{n-3}+\cdots]\nonumber\\
=& yx^{n-3} +\cdots\label{eq-B7}
\end{align}
since
\begin{align*}
&F(C_{k_1, k_2, k_3}; \omega_2)  \\
=& (-1)^{k_1}\Big[\frac{(1-\omega_2)^{k_1}-1}{\omega_2}F(D_{k_2+k_3};\omega_2) + F(P_{k_2, k_3};\omega_2)\Big]\\
=& (-1)^{k_1}[(-1)\sum_{i=0}^{k_1-1}y^i(-1)^{k_2+k_3-1}\sum_{i=1}^{k_2+k_3-1}y^{i} + (-1)^{k_2+k_3}\sum_{i=1}^{k_2-1}y^{i}\sum_{i=1}^{k_3-1}y^{i}]\\
=& (-1)^{k}(y + 3y^2 + \cdots ). \\
\end{align*}

Now let's insert (\ref{eq-B1})-(\ref{eq-B7}) into (\ref{eq-sum}) and will obtain

\begin{align*}
&T(G; x, y)\nonumber\\
=&x^{n-1} + (-1)[(n-1) + \sum_{\stackrel {H\in L(G)}{rk(H)=1}}\mu(\widehat{0},H)]x^{n-2}  \\
& + [{n-1 \choose 2} + (n-2)\sum_{\stackrel {H\in L(G)}{rk(H)=1}}\mu(\widehat{0},H) + \sum_{\stackrel {H\in L(G)}{rk(H)=2}}\mu(\widehat{0},H)]x^{n-3}\\
&+ \cdots \\
&+ \sum_{\stackrel {A \subseteq E }{A\ \text{is a nontrivial parallel class}}} yx^{n-2}\\
& + \Big\{\sum_{\stackrel {A \subseteq E }{A\  \text{is a nontrivial parallel class}}}[-(n-2)-\sum_{\stackrel {H\in L(G/A)}{rk(H)=1}}\mu(\widehat{0},H)]+ \sum_{\stackrel {A \subseteq E,G[A]=C_{k_1,k_2,k_3}} {G/A\ \text{is loopless}}}1 \Big\}yx^{n-3}\\
& + \cdots.
\end{align*}

Clearly,
\begin{eqnarray*}
\sum_{\stackrel {H\in L(G)}{rk(H)=1}} \mu(\widehat{0}, H)&=&-p(G),\\
\sum_{\stackrel {A \subseteq E,G[A]=C_{k_1,k_2,k_3}} {G/A\ \text{is loopless}}}1&=&\Delta(\tilde{G}).
\end{eqnarray*}

Let $\tilde{G}$ be the graph obtained from $G$ by replacing each parallel class by a single edge. Note that $rk(H)=2$ means $c(H)=n-2$. $H=P_3\cup E_{n-3}$ (but $\tilde{G}[V(P_3)]\neq C_3$), $P_2\cup P_2\cup E_{n-4}$ or $C_3\cup E_{n-3}$. The former two have  the M\"{o}bius function value 1 and the third one has  the M\"{o}bius function value 2. Then
\begin{eqnarray*}
\sum_{\stackrel {H\in L(G)}{rk(H)=2}} \mu(\widehat{0}, H)=1\times[{p(G) \choose 2}-3\Delta(\tilde{G})]+2\times\Delta(\tilde{G})={p(G) \choose 2}-\Delta(\tilde{G}).
\end{eqnarray*}
This completes the proof of Theorem \ref{mainII}.
\hfill\cvd

\section{Discussions}
\noindent

In this section, we first discuss the duality of Theorems 1 and 2. Then, we deduce the results on extreme coefficients of the Jones polynomials of alternating links \cite{Lin} and graphs \cite{DJIN}. Finally we also deduce extreme coefficients of the chromatic and flow polynomials.

\subsection{Duality}
\noindent

It is well known that  if $G$ is a plane graph and $G^*$ is the dual graph of $G$, then
\begin{eqnarray}
T(G^*;x,y)=T(G;y,x).
\end{eqnarray}

Let $t^*_{i,j}$ be the coefficient of $x^iy^j$ in the Tutte polynomial $T(G^*;x,y)$, and we have $t^*_{i,j}=t_{j,i}$. Loops and bridges, deletion and contraction will interchange by taking dual of a plane graph. The dual of a bridgeless and loopless connected plane graph is still a bridgeless and loopless connected plane graph. Let $G$ be a bridgeless and loopless connected plane graph and $G^*$ be the dual of $G$. Let $n^*$ and $m^*$ be the order and size of $G^*$, respectively. Then

\begin{eqnarray*}
m^*&=&m,\\
n^*&=&m-n+2,\\
s(G^*)&=&p(G),\\
p(G^*)&=&s(G),\\
s^*(G^*)&=&p^*(G),\\
p^*(G^*)&=&s^*(G),\\
\Theta(G^*)&=&\Delta(\tilde{G}),\\
\Delta(\tilde{G^*})&=&\Theta(G).
\end{eqnarray*}

Theorems \ref{mainI} and \ref{mainII} are dual to each other in the case that $G$ is a bridgeless and loopless connected plane graph. Now we check it as follows.

\begin{itemize}
\item[(1)]
\begin{eqnarray*}
t^*_{0,n-1}=t^*_{0,m^*-n^*+1}=1=t_{n-1,0}.
\end{eqnarray*}
\item[(2)]
\begin{eqnarray*}
t^*_{0,n-2}&=&t^*_{0,m^*-n^*}\\
&=&n^*+s(G^*)-m^*-1\\
&=&(m-n+2)+p(G)-m-1\\
&=&p(G)-n+1\\
&=&t_{n-2,0}.
\end{eqnarray*}
\item[(3)]
\begin{eqnarray*}
t^*_{0,n-3}&=&t^*_{0,m^*-n^*-1}\\
&=&{m^*-n^*+1 \choose 2} - (m^*-n^*)s(G^*) + {s(G^*) \choose 2} - \Theta(G^*)\\
&=&{n-1 \choose 2} - (n-2)p(G) + {p(G) \choose 2} - \Delta(\tilde{G})\\
&=&t_{n-3,0}.
\end{eqnarray*}
\item[(4)]
\begin{eqnarray*}
t^*_{1,n-2}&=&t^*_{1,m^*-n^*}\\
&=&s^*(G^*)\\
&=&p^*(G)\\
&=&t_{n-2,1}.
\end{eqnarray*}
\item[(5)]
\begin{eqnarray*}
t^*_{1,n-3}&=&t^*_{1,m^*-n^*-1}\\
&=&-s^*(G^*)(m^*-n^*)+\sum_{\stackrel{A^*\subseteq E^*}{A^*\ \text{is a nontrivial series class} }  }s(G^*-{A^*}) + \Theta(G^*)\\
&=& -p^*(G)(n-2)+\sum_{\stackrel{A\subseteq E}{ A\ \text{is a nontrivial  parallel class} }  }p(G/A) + \Delta(\tilde{G})\\
&=&t_{n-3,1}.
\end{eqnarray*}
\end{itemize}

In addition, Theorems 1 and 2 may be generalized to the Tutte polynomials of matroids and in that setting the duality may be more obvious.

\subsection{Jones polynomial}
\noindent

In \cite{DJIN}, Dong and Jin introduced the Jones polynomial of graphs. In the case of plane graphs, it (up to a factor) reduces to the Jones polynomial of the alternating link constructed from the plane graph via medial construction \cite{Th,BOLLOBAS}.
We denote by $J_G(t)$ the Jones polynomial of $G$. Then
\begin{eqnarray}
J_G(t)=(-1)^{n-1}t^{m-n+1}T(-t,-t^{-1}).
\end{eqnarray}

Based on the work of Dasbach and Lin \cite{Lin}, Dong and Jin \cite{DJIN} further obtain:

\begin{theorem}[\cite{DJIN}]\label{cor}
Let $G=(V,E)$ be a connected bridgeless and loopless graph with order $n$ and size $m$. Then \[J_G(t)=b_0+b_1t+b_2t^2+...+b_{m-2}t^{m-2}+b_{m-1}t^{m-1}+b_{m}t^{m},\]
where $(-1)^{m-i}b_i$ is a non-negative integer for $i=0,1,2,\cdots,m$ and in particular,
\begin{eqnarray*}
b_0&=&(-1)^{m},\\
b_1&=&(-1)^{m}[m-n+1-s(G)],\\
b_{m-2}&=&\binom{p(G)-n+2}{2}+p^*(G)-\Delta(\tilde{G}),\\
b_{m-1}&=&n-1-p(G),\\
b_{m}&=&1.
\end{eqnarray*}
\end{theorem}

We can deduce Theorem \ref{cor} by using Theorem \ref{mainI} and Theorem \ref{mainII} and taking $x=-t,y=-t^{-1}$, and further obtain:
\begin{corollary}
\begin{eqnarray*}
b_2=(-1)^{m}\left[\binom{s(G)-m+n}{2}+s^*(G)-\Theta(G)\right].
\end{eqnarray*}
\end{corollary}

\noindent\textbf{Proof.}
\begin{eqnarray*}
b_2&=&(-1)^{m}[t_{0,m-n-1}+t_{1,m-n}]\\
&=&(-1)^{m}\left[{m-n+1 \choose 2} - (m-n)s(G) + {s(G) \choose 2} - \Theta(G)+s^*(G)\right]\\
&=&(-1)^{m}\left[\binom{s(G)-m+n}{2}+s^*(G)- \Theta(G)\right].
\end{eqnarray*} \hfill\cvd

In \cite{Kau}, Kauffman generalized the Tutte polynomials from graphs to signed graphs, which includes the Jones polynomial of both alternating and non-alternating links. It is worth studying extreme coefficients of the signed Tutte polynomial.

\subsection{Chromatic and flow coefficients}
\noindent

\begin{theorem}[\cite{Read,Mer}]
Let $G$ be a loopless and bridgeless connected graph of order $n$. Let $P(G;\lambda)=a_0\lambda^n+a_1\lambda^{n-1}+\cdots+a_{n-1}\lambda$. Then
\begin{eqnarray*}
a_0&=&1,\\
a_1&=&-p(G),\\
a_2&=&\binom{p(G)}{2}-\Delta(\tilde{G}).
\end{eqnarray*}
\end{theorem}

By Eq. (\ref{ccc}), one can obtain $t_{n-1,0}, t_{n-2,0}$ and $t_{n-3,0}$ and vice versa. Theorem \ref{fm} may be known but we have not found it in the literature. By Eq. (\ref{eq-A1}), one can obtain it from $t_{0,m-n+1}, t_{0,m-n}$ and $t_{0,m-n-1}$.

\begin{theorem}\label{fm}
Let $G$ be a loopless and bridgeless connected graph of order $n$ and size $m$. Let $F(G;\lambda)=c_0\lambda^{m-n+1}+c_1\lambda^{m-n}+\cdots+c_{m-n+1}$. Then
\begin{eqnarray*}
c_0&=&1,\\
c_1&=&-s(G),\\
c_2&=&\binom{s(G)}{2}-\Theta(G).
\end{eqnarray*}
\end{theorem}

\section*{Acknowledgements}
\noindent

This work is supported by NSFC (No. 11271307,11671336) and President's Funds of Xiamen University (No. 20720160011). We thank the anonymous referee and A/P Fengming Dong for some helpful comments.

\end{document}